\documentclass[12pt,a4paper]{amsart}
\usepackage[utf8]{inputenc}	
\usepackage[T1]{fontenc}
\usepackage[in,headings]{fullpage}
\usepackage{amsbsy, amscd, amsfonts, amsmath, amsrefs, amssymb, amsthm}
\usepackage{graphicx}
\usepackage{float}
\usepackage{color}   
\usepackage[colorlinks=true,linkcolor=blue,citecolor=blue,urlcolor=blue]{hyperref}

\newtheorem{theorem}{Theorem}

\newtheorem{lemma}[theorem]{Lemma}

\newtheorem*{conjecture}{Conjecture}
\newtheorem*{question}{Question}
\newtheorem{problem}{Problem}

\theoremstyle{definition}

\theoremstyle{remark}

\newcommand{\C}{\mathbf{C}}

\newcommand{\R}{\mathbf{R}}

\renewcommand{\Re}{\mathop{\mathrm{Re}}\nolimits}
\renewcommand{\Im}{\mathop{\mathrm{Im}}\nolimits}
\newcommand{\Rzeta}{\mathop{\mathcal R }\nolimits}
\newcommand{\Lzeta}{\mathop{\mathcal L }\nolimits}
\newfont{\cmbsy}{cmbsy10}
\newfont{\cmmib}{cmmib10}
\newcommand{\Orden}{\mathop{\hbox{\cmbsy O}}\nolimits}
\newcommand{\orden}{\mathop{\hbox{\cmmib o}}\nolimits}

\begin{document}

\title[Papers related to Riemann's auxiliary function]%
{Report on some papers related to the function $\Rzeta(s)$ found by Siegel in Riemann's posthumous papers.}
\author[Arias de Reyna]{J.~Arias de Reyna}
\address{%
Universidad de Sevilla \\ 
Facultad de Matem\'aticas \\ 
c/Tarfia, sn \\ 
41012-Sevilla \\ 
Spain.} 

\subjclass[2020]{Primary 11M06; Secondary 30D99}

\keywords{zeta function, Riemann's auxiliar function}


\email{arias@us.es, ariasdereyna1947@gmail.com}


\begin{abstract}
In a letter to Weierstrass Riemann asserted that the number $N_0(T)$ of zeros of $\zeta(s)$ on the critical line to height $T$ is approximately equal to the total number of zeros to this height $N(T)$. Siegel studied some posthumous papers of Riemann trying to find a proof of this. He found a function $\Rzeta(s)$ whose zeros are related to the zeros of the function $\zeta(s)$. Siegel concluded that Riemann's papers contained no ideas for a proof of his assertion, connected the position of the zeros of $\Rzeta(s)$ with the position of the zeros of $\zeta(s)$ and asked about the position of the zeros of $\Rzeta(s)$. This paper is a summary of several papers that we will soon upload to arXiv, in which we try to answer Siegel's question about the position of the zeros of $\Rzeta(s)$. The articles contain also improvements on Siegel's results and also other possible ways to prove Riemann's assertion, but without achieving this goal.
\end{abstract}

\maketitle

\section{Introduction}
Let $N(T)$ be the number of zeros of $\zeta(s)$ in the rectangle $(0,1)\times(0,T]$ counted with multiplicity and $N_0(T)$ the number of those zeros $\varrho$ with $0<\Im\varrho\le T$ and on the critical line. Riemann claims in a letter to Weierstrass that he has a  relatively complicated proof that $N_0(T)$ is approximately equal to $N(T)$.  Therefore, we state here the conjecture.
\begin{conjecture}[Riemann's last Theorem] The number of zeros $N_0(T)$ of $\zeta(s)$ on the critical line with $0\le \Im(\rho)\le T$ is $N_0(T)=\frac{T}{2\pi}\log\frac{T}{2\pi}-\frac{T}{2\pi}+\orden(T)$ for $T\to+\infty$.
\end{conjecture}

Riemann phrase in his paper \cite{R}: \emph{One finds in fact about this many real roots between these bounds and it is very likely that all of the roots are real},  is sufficiently vague for us to think that he means that he has calculated the first three or four zeros and shown them to be real. Annotations of these calculations have been found in his posthumous papers. In a letter to Weierstrass \cite{R2}*{p.~823} he is much more explicit.

\begin{quote}
I have not yet fully carried out the proof; and about this proof 
I would like to note --- you may allow me to refer to the papers
to be sent to you very soon ---  that
the two propositions that I only have announced:
\medskip

\noindent \emph{That between $0$ and $T$ there are approximately $\frac{T}{2\pi}\log\frac{T}{2\pi}-\frac{T}{2\pi}$ real roots of the equation $\xi(\alpha)=0$. }
\medskip

and [another one\dots]
 
\noindent follow from a new development of the function $\xi$, which I have not yet sufficiently simplified to announce it. From the correctness of everything else, although some are  only very briefly indicated, you will easily convince yourself.
\end{quote}

An expansion of the function $\xi$ was found in the papers of Riemann, and Siegel \cite{Siegel} tried to see what he could get out of this expression. He found the Riemann-Siegel expansion and studied the function  $\Rzeta(s)$ in connection with the zeros of zeta. By means of $\Rzeta(s)$ he obtained that $N_0(T)\ge CT$ for some explicit constant $C$, but concluded that there was nothing in Riemann's papers in the direction of proving the above conjecture. He asked about the position of the zeros of $\Rzeta(s)$, and some of my work addresses this question of Siegel.

Riemann's papers on Number Theory in Göttingen are not the papers that a mathematician would keep from his research; rather, they are the ones he would throw in the wastebasket after writing down the most important things in his notes. If Riemann had any ideas for proving his conjecture, they are not contained in those papers (as Siegel concluded in \cite{Siegel}).  What is clear is that Riemann's papers show that he knew much more about the zeta function than is contained in his article \cite{R}.

Siegel introduced the function
\begin{equation}\label{R:def}
\Rzeta(s)=\int_{0\swarrow1}\frac{x^{-s} e^{\pi i x^2}}{e^{\pi i x}-
e^{-\pi i x}}\,dx.
\end{equation}
We will call this Riemann's auxiliary function. The papers we present study this function, refining Siegel's results in \cite{Siegel}, trying to resolve Siegel's question about the location of its zeros, and asking about some issues that we have not been able to resolve.  

Apart from this, some of the papers we present try to explore other ways of proving the Riemann conjecture. In this direction we must cite the work of Levinson \cite{Lev} who proved $N_0(T)\ge\frac13N(T)$. And the improvements by Conrey \cite{Con}, and several others, for example Pratt, Robles, Zaharescu and Zeindler \cite{Pr} prove that $N_0(T)\ge\frac{5}{12}N(T)$ and speculate about the possibility to get Riemann's last theorem by refining these techniques.

Siegel (or Riemann) related $\Rzeta(s)$ to the zeta function by means of the equation
\begin{equation}\label{zetaRzeta}
\zeta(s)=\Rzeta(s)+\chi(s)\overline{\Rzeta}(1-s);\qquad
\chi(s)=\pi^{s-1/2}
\frac{\Gamma\left(\frac{1-s}{2}\right)}{\Gamma\left(\frac{s}{2}\right)},\quad
\overline{\Rzeta}(s)=\overline{\Rzeta(\overline{s})}.
\end{equation}
In particular, we have the following expression for the Riemann-Siegel function $Z(t)$\footnote{I am reluctant to call this Hardy's function, after having seen it in Riemann's papers at Göttingen, and in his commentary by Siegel.}
\begin{equation}
Z(t)=2\Re\{e^{i\vartheta(t)}\Rzeta(\tfrac12+it)\},\text{ where } \zeta(\tfrac12+it)=e^{-i\vartheta(t)}Z(t)\qquad t\in\R.
\end{equation}
 
The following is a summary of each of the papers, containing our results on $\Rzeta(s)$.

\section{Riemann-Siegel formula}

\subsection*{1. High precision computation of Riemann's zeta function by the Riemann-Siegel formula I \texorpdfstring{\cite{A86}}{[3]}} Siegel obtained the Riemann-Siegel asymptotic expansion without giving explicit error bounds, Gabcke in his thesis \cite{Gabcke} obtained explicit bounds for $\sigma=\frac12$ and after retaining less than ten terms. In this paper, we follow the steps of Gabcke to obtain realistic explicit bounds after retaining any number of terms at almost any point in the plane. 

In \cite{A86}*{Prop.~6.1} the author needs to prove the inequality
\begin{equation}\label{E:RS-ineq}
\sum_{k=1}^\infty(-1)^{k+1}\frac{r^k\cos k\theta}{k+2}\le \sum_{k=1}^\infty(-1)^{k+1}\frac{r^k}{k+2},\end{equation}
for $r=1$, $r=8/9$ and $r=0.883$, and prove it by means of the maximal slope principle. After \cite{A86} was published, J. van de Lune and J. Arias de Reyna \cite{A94} proved this inequality for any $r$. Later on, Gabcke \cite{G2} gave a nice proof of the inequality \eqref{E:RS-ineq}. 

\subsection*{2. High precision computation of Riemann's zeta function by the Riemann-Siegel formula II \texorpdfstring{\cite{A88}}{[4]}} After writing \cite{A86} we wanted to implement in Python the computation of $\Rzeta(s)$ and $\zeta(s)$. Given $s\ne 1$ and $\varepsilon$, the problem was to compute $z\in\C$ such that $|\zeta(s)-z|<\varepsilon$. Is this possible? How many terms of the expansion have to compute? To what precision? We solved this thanks to the power of \texttt{mpmath} a library for arbitrary-precision arithmetic created by F. Johansson \cite{mpmath}. In \texttt{mpmath} my program can be used freely. This paper solves these questions. But it is not written in final form. 

We worked from scratch. We noticed that the usual implementations of the product of two complex numbers do not follow the IEEE standard. Given two complex numbers $a$, $b$, we want the approximation $c$ given by the product to satisfy $|c-ab|\le |ab| 2^{-d}$ when working with $d$ binary digits. Our calculations would be simpler if this rule were followed.

\subsection*{3. Programs for Riemann's zeta function \texorpdfstring{\cite{A116}}{[8]}} This article deals with the extension of programs to compute $\Rzeta(s)$ and $\zeta(s)$ in \texttt{mpmath} \cite{mpmath}, to the computation of the first four derivatives of zeta and the $n$-th zero given $n$. Computing $\varrho_n$ starting from $n$ is not as easy as some pretend, the greatest difficulty lies in precisely isolating the $n$-th zero and not any other.

\section{Riemann auxiliary function generalities}

\subsection*{4. Riemann auxiliary function. Basic results \texorpdfstring{\cite{A166}}{[9]}} This paper contains the definition \eqref{R:def} and some other representation, such as, for example,
\begin{equation}\label{E:Rzetatres}
\pi^{-s/2}\Gamma(s/2)\Rzeta(s)=-\frac{e^{-\pi i s/4}}{ s}
\int_{-1}^{-1+i\infty} \tau^{s/2}\vartheta_3'(\tau)\,d\tau.
\end{equation}
In addition, we determine the values of $\Rzeta(n)$ for any even positive number and for odd negative integers. The trivial zeros of zeta are also zeros of $\Rzeta$. 

\begin{figure}[H]
\begin{center}
\includegraphics[width=\hsize]{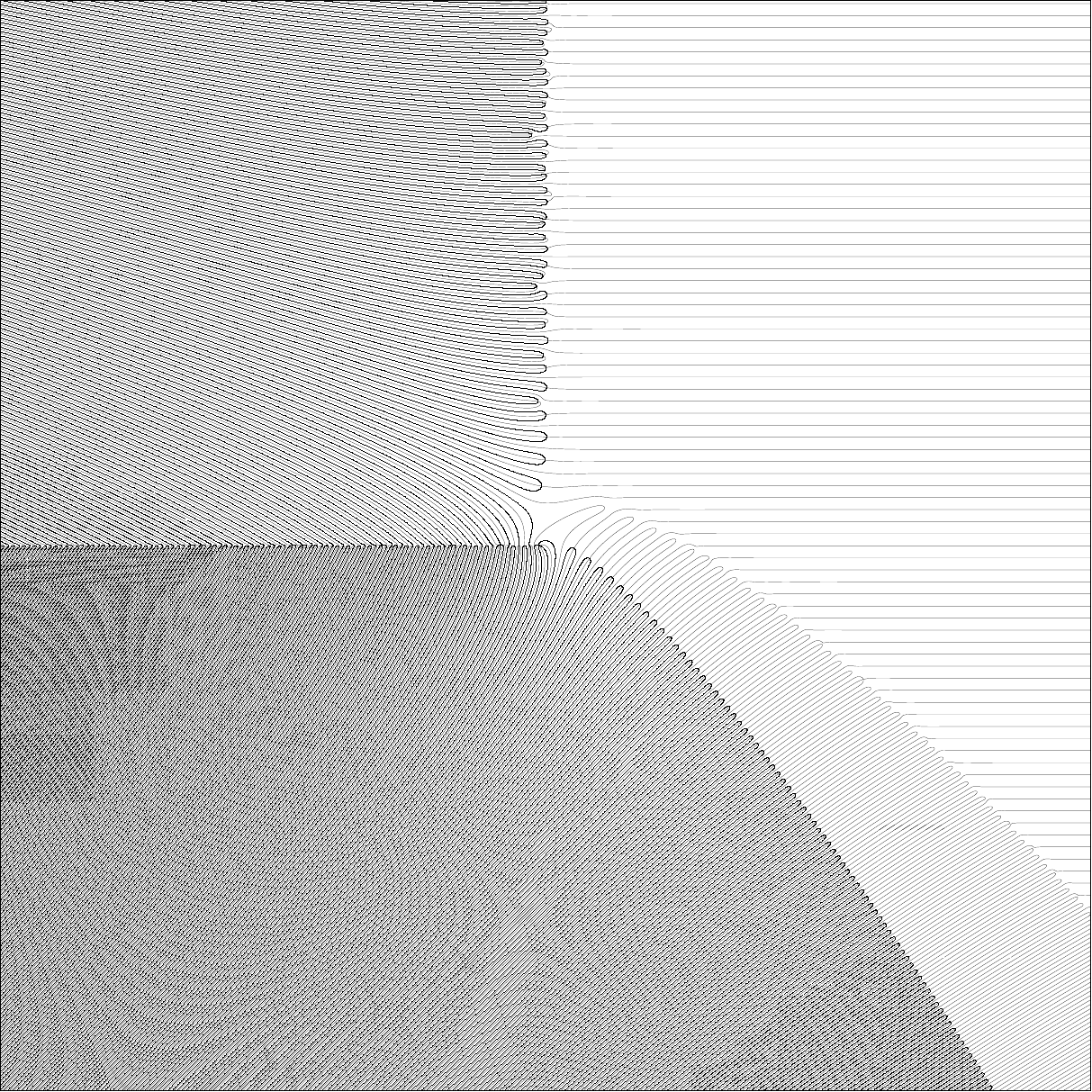}
\caption{x-ray of $\Rzeta(s)$ on $(-200,200)^2$ }
\label{rzeta}
\end{center}
\end{figure}

Siegel wondered where the zeros of the function $\Rzeta(s)$ were. Some of my works will provide a partial answer to this question. The X-ray of a function consists of drawing the lines of the plane where the function takes real values, or purely imaginary values. Figure \ref{rzeta} shows the X-ray\footnote{In the pdf an X-ray can be magnificated to see the details.} of function $\Rzeta(s)$ in rectangle $(-200,200)^2$. The computation of $\Rzeta(s)$ for moderate values of $s$ is difficult. To obtain part of this X-ray we have needed to calculate $\Rzeta(s)$ with 80  decimal digits.  This X-ray shows the distribution of the zeros of $\Rzeta(s)$ into three main types. More on this in some of the papers.

\subsection*{5. Statistics of zeros of Riemann's auxiliary function \texorpdfstring{\cite{A172}}{[10]}}
We have computed all zeros $\beta+i\gamma$ of $\Rzeta(s)$ with $0<\gamma<215946.3$. A total of 162215 zeros with 25 correct decimal digits. The zeros with $\gamma>0$ were the most interesting zeros for Siegel. $63.9\%$ of these zeros satisfies $\beta<\frac12$. If this were true for all zeros of $\Rzeta(s)$ with $\gamma>0$, it will follow that $63.9\%$ of zeros of $\zeta(s)$ will be on the critical line. Therefore, it is plausible that the percentage of zeta zeros on the critical line can be improved by studying $\Rzeta(s)$. But this author has not been able to do that.

There are other interesting things in this paper. For example, studying these zeros, we may guess  that the number of zeros of $\Rzeta(s)$ with $0<\gamma\le T$ satisfies
\begin{equation}\label{E:N(T)}
N(T)=\frac{T}{4\pi}\log\frac{T}{2\pi}-\frac{T}{4\pi}-\frac12\sqrt{\frac{T}{2\pi}}+\frac32+O(\log T).
\end{equation}
In \cite{A185} we prove this by substituting the error $\frac32+O(\log T)$ by $\Orden(T^{2/5}\log^2T)$. We also observe a periodic structure on the zeros with $\gamma>0$. 

Another interesting observation was that zeros with $\gamma>0$ satisfies $\beta<1$. More on this in \cite{A173}.
\vfill

\section{Asymptotic expansion for the auxiliary function}

\subsection*{6. Region without zeros for the auxiliary function of Riemann \texorpdfstring{\cite{A98}}{[11]}} This paper focuses on the  Siegel result \cite{Siegel} about the left limit of the zeros of $\Rzeta(s)$ for $t>0$. We extend the  asymptotic expansion of Siegel in the second quadrant, extending them to most of the third quadrant. Giving precise bounds of the error allows us to give an explicit region free of zeros or with only trivial zeros. 

There is an interesting question here. Siegel claims that given $\varepsilon>0$ there is a $t_0>0$ such that there is no zero of $\Rzeta(s)$ with  $1-\sigma> ct^\varepsilon$ and $t>t_0$. Siegel only proves the case $\varepsilon=3/7$. We were unable to extend Siegel's reasoning beyond $1-\sigma>t^{2/5}\log t$. This limits the error, for example, in \eqref{E:N(T)}. Therefore,  my question:
\begin{question} \label{Question}
Is it true that the zeros $\beta+i\gamma$ of $\Rzeta$ with $\gamma>0$ satisfy  $\limsup_{\gamma\to+\infty}\frac{1-\beta}{\gamma^\varepsilon}=0$ for any $\varepsilon>0$?
\end{question}

\subsection*{7.~Asymptotic expansions of the auxiliary function \texorpdfstring{\cite{A100}}{[12]}} The terms of the asymptotic expansion in \cite{A86} are not analytic functions. This simplifies  obtaining the expansion of the derivatives of $\Rzeta(s)$ (and consequently those of $\zeta(s)$ or $Z(t)$). We have used this in our programs in mpmath, but have not written a paper explaining this. By comparison, the terms of the expansion obtained by Siegel are analytic functions. Here we give two related asymptotic expansions for $\Rzeta(s)$; this is only suggested in Siegel, and we also extend the range of validity, specifying regions without zeros. 
For example, we get for $t\to+\infty$
\[e^{-i\vartheta(t)}\Rzeta(\tfrac12-it)=
-\frac{1}{\sqrt{2}}\Bigl(\frac{t}{2\pi}\Bigr)^{-\frac14}
\exp\Bigl\{\frac{\pi t}{2}-\Bigl(\frac{\pi t}{2}\Bigr)^{\frac12}\Bigr\}\cdot i\cdot(1+\Orden(t^{-1/2})),\]
this is almost purely imaginary and very large. The real part is $Z(t)$, which is relatively small. Therefore, for $t<0$ our expression for $Z(t)$ is almost useless to locate zeros. 

In this paper, we show that $\Rzeta(s)$ is approximated by an adequate zeta sum in the first and part of the fourth quadrant (see figure \ref{rzeta}). There is a narrow zone in the fourth quadrant limit between the two asymptotic expansions where a line of zeros of $\Rzeta(s)$ is found (again, see Figure \ref{rzeta}). Finally, we explain the somewhat inexcrutable formulas (84)  and (85) of Siegel in \cite{Siegel}. 

\subsection*{8.~Note on the asymptotic part of the auxiliary function \texorpdfstring{\cite{A193}}{[14]}} To define explicit regions without zeros in a first paper, we get an approximation of $\Rzeta(s)$ of type $f(s)(1+U)$ with $|U|< 1$. But this $U$ does not tend to zero when $t\to+\infty$. In the present paper we get an approximation of the form $f(s)(1+\orden(t))$. We clarify here Siegel's result following his reasoning. This is essential to get the last Theorems in Siegel's paper about $\Rzeta(s)$. We prove 

\begin{theorem}
There exist constants $A$ and $t_0>1$ such that for $s$ in the closed set
\[\Omega=\{s\in\C\colon t\ge t_0,\quad 1-\sigma\ge t^{3/7}\},\]
we have 
\begin{equation}\label{E:withUT3}
\Rzeta(s)=-\chi(s)\eta^{s-1}e^{-\pi i \eta^2}\frac{\sqrt{2}e^{3\pi i/8}\sin\pi\eta}{2\cos2\pi\eta}(1+U),\qquad |U|\le At^{-\frac{1}{21}}.
\end{equation}
\end{theorem}

\section{Zeros of \texorpdfstring{$\Rzeta(s)$}{R(s)}}

\subsection*{9.~Riemann auxiliary function. Right limit of zeros \texorpdfstring{\cite{A173}}{[13]}}
All the zeros $\beta+i\gamma$ of $\Rzeta(s)$ that we computed satisfies $\beta<1$.  Here, we prove that this is true for $\gamma>t_0=3.93\times10^{65}$. We use the paper \cite{A86} to get a good approximation to $\Rzeta(s)$ by zeta sums. This can be a model of how to use \cite{A86}. Of course, it will be desirable to prove this for $t_0=0$, but this appears to be difficult. 

We need explicit bounds for zeta sums. It has been unfortunate that a resuly by Landau from 1928 has been widely overlooked, rendering a number of less-than-optimal results. We clarify this point in the paper \cite{A182}, and then obtained explicit results for the van der Corput $d$-th derivative test. This is the origin of my work in \cite{A181}.

\subsection*{10.~On Siegel results about the zeros of the auxiliary function of Riemann \texorpdfstring{\cite{A102}}{[16]}} In this paper, we fill the details on the reasoning at the end of Siegel's paper \cite{Siegel}. But we feel free to change some things when we can simplify the reasoning. At some points, we went a little further.
For example, Siegel proved  
\[\sum_{0<\gamma\le T}\beta=-\frac{T}{4\pi}\log2+\orden(T).\]
we improve this in \cite{A102}*{Prop.~6} giving the error $\Orden(T^{20/21})$. But we suspect that if the Siegel claim about the zeros (Question \ref{Question}) is true, we can get an additional term that will explain the experimental results in \cite{A172}.

\subsection*{11.~On the number of zeros of \texorpdfstring{$\Rzeta(s)$}{R(s)} \texorpdfstring{\cite{A185}}{[17]}} We prove the formula \eqref{E:N(T)} for the number of zeros $N(T)$ for $\Rzeta(s)$ with $0<\gamma\le T$. This improves Siegel's result and explains the experimental results seen in \cite{A172}.

I take this opportunity to give a reasonable proof of the Backlund Lemma \cite{B}
\begin{lemma}[Backlund]
Let  $f$ be holomorphic in the disc $|z-a|\le R$. 
Let  $|f(z)|\le M$ for  $|z-a|\le R$. Let  $b$ be a point in the interior of the disc
$0<|b-a|<R$. Assume that  $f$ does not vanish on the segment $[a,b]$, then 
\begin{equation}\label{E:backlund}
\Bigl|\Re\frac{1}{2\pi i}\int_a^b\frac{f'(z)}{f(z)}\,dz\Bigr|\le
\frac{1}{2}\Bigl(\log \frac{M}{|f(a)|}\Bigr)\Bigl(\log\frac{R}{|b-a|}\Bigr)^{-1}.
\end{equation}
\end{lemma}

\subsection*{12.~Trivial zeros of Riemann auxiliary function \texorpdfstring{\cite{A186}}{[18]}} In this paper, it is proved that the trivial zeros of $\Rzeta(s)$ are simple. We get an integral expression for the derivative
\begin{equation}\label{E:main}
\Rzeta'(-2n)=\frac{\omega\sqrt{\pi n}}{2}\Bigl(\frac{i n}{\pi e}\Bigr)^n\int_0^\infty\frac{(ex^2e^{-x^2})^n}{\sin(\omega x\sqrt{\pi n})}\,dx,
\end{equation}
and show that these numbers are not $0$.

\subsection*{13.~Zeros of \texorpdfstring{$\Rzeta(s)$}{R(s)} in the fourth quadrant \texorpdfstring{\cite{A108}}{[19]}} We show that there exists a sequence $\rho_{-1}$, $\rho_{-2}$, \dots\ of zeros of $\Rzeta(s)$ in the fourth quadrant. We give an algorithm to compute the $n$-th zero on the fourth quadrant, and an asymptotic expansion for it showing that
\[\rho_{-n}\sim \frac{4\pi^2n}{\log^2n}-\frac{4\pi i n}{\log n}.\]
This is needed in \cite{A66}.

\subsection*{14.~Infinite product of Riemann auxiliary function \texorpdfstring{\cite{A66}}{[20]}} We obtain the product formula for $\Rzeta(s)$. We study the phase of $\Rzeta(s)$ at the critical line, that is, the real analytic  and real function $\omega(t)$ such that $\Rzeta(\frac12+it)=e^{-i\omega(t)}g(t)$ where $g$ is also real and real analytic. Each zero of $\cos(\vartheta(t)-\omega(t))$ determine a zero of $\zeta(s)$ on the critical line. We study the relation of $\omega(t)$ with the zeros of $\Rzeta(s)$. For example, for $t>0$ we have $\omega(t)\approx2\pi N_r(t)$ where $N_r(t)$ is the number of zeros $\rho=\beta+i\gamma$ with $0<\gamma<T$ and $\beta>\frac12$. 

Our computation of zeros (see \cite{A172}) of $\Rzeta(s)$  suggests that $N_r(t)\approx \frac{\vartheta(t)}{6\pi}$. If this were true $\vartheta(t)-\omega(t)\approx 
\frac{2}{3}\vartheta(t)$, and at least $\frac23$ of the zeros of $\zeta(s)$ will be on the critical line. The method of Levinson \cite{Lev}, has been refined by several authors  \cite{Con}, \cite{F}, \cite{Pr}, to obtain $\frac{5}{12}$ zeros of $\zeta(s)$ on the critical line. Therefore, the study of the zeros of $\Rzeta(s)$  has the potential  to improve these results. 

The behavior of $\Rzeta(s)$ on the critical line, similar to that of $\zeta(s)$, depends on the sum \[\sum_{n\le\sqrt{t/2\pi}} n^{-s}.\] We know that the behavior of $\zeta(s)$ is not typical until the $\log\log t$ is really large. Therefore, assume that $N_r(t)\approx \frac{\vartheta(t)}{6\pi}$ from the scarce data we have is risky. Maybe $N_r(t)=\orden(t) $ after all, and Riemann had a proof, who knows.

\section{Density Theorems}

\subsection*{15. Mean Values of the auxiliary function \texorpdfstring{\cite{A101}}{[15]}}
In this paper, we obtain the main terms of the mean values 
\[\frac{1}{T}\int_0^T |\Rzeta(\sigma+it)|^2\Bigl(\frac{t}{2\pi}\Bigr)^\sigma\,dt, \quad\text{and}\quad \frac{1}{T}\int_0^T |\Rzeta(\sigma+it)|^2\,dt.\]
Giving complete proofs of some result of the paper of Siegel on the Riemann Nachlass.  Here, $\Rzeta(s)$ is the function related to $\zeta(s)$ found by Siegel in the Riemann papers.

There are two cryptic passages in Siegel \cite{Siegel}*{before eq. (64)} where he says 
\begin{quote}
Riemann tried to obtain an expression about the zeros of $\Rzeta(s)$; from (58) he found 
\begin{equation}\label{E:doubleInt}
|\Rzeta(\sigma +ti)|^2=\int\limits_{0\swarrow1}\int\limits_{0\searrow1}
\frac{x^{-\sigma-ti}y^{-\sigma+ti}e^{\pi i(x^2-y^2)}}{ (e^{\pi i
x}-e^{-\pi i x})(e^{\pi i y}-e^{-\pi i y})}\,dx\,dy,
\end{equation}
and brought the complex double integral into a different form by the introduction of new variables, the deformation of the integration area, and the use of the residue theorem; however, this did not yield a useful result.
\end{quote}
then in \cite{Siegel}*{between eq. (86) and (87)} to compute a mean value he writes
\begin{quote}
\[\int_0^\infty |\Rzeta(\sigma+it)|^2e^{-\varepsilon t}\,dt=
\int_0^\infty e^{-\varepsilon t}\Bigl\{\int\limits_{0\swarrow1}\int\limits_{0\searrow1}
\frac{x^{-\sigma-ti}y^{-\sigma+ti}e^{\pi i(x^2-y^2)}}{ (e^{\pi i
x}-e^{-\pi i x})(e^{\pi i y}-e^{-\pi i y})}\,dx\,dy\Bigr\}\,dt.\]
Here, one may transform the right side by deformation of the path of integration, interchange of the order of integration, and application of the residue theorem. The calculation yields the statement
\begin{equation}\label{firstequation} 
\int_0^\infty |\Rzeta(\sigma+it)|^2e^{-\varepsilon t}\,dt\sim\frac{1}{2\varepsilon}(2\pi\varepsilon)^{\sigma-\frac12}\Gamma(\tfrac12-\sigma).
\end{equation}
\end{quote}
We were unable to reproduce any of this. So, we get the result by more standard means.

From \eqref{firstequation} Siegel asserts that for all
$\sigma<\frac12$
\begin{equation}\label{forsigma}
\int_1^{\infty}|\Rzeta(\sigma+it)|^2\Bigl(\frac{t}{2\pi}\Bigr)^\sigma 
e^{-\varepsilon t}\,dt\sim \frac{(2\varepsilon)^{-\frac32}}{1-2\sigma}.
\end{equation}
He does not give a proof of this implication. There is a statement in  
Titchmarsh \cite{T}*{p.~159} which proves the implication \eqref{firstequation} $\Longrightarrow$ \eqref{forsigma} for $\sigma<0$.  The proof in Titchmarsh needs  hypothesis $\sigma<0$.  Siegel only uses \eqref{forsigma} for a value of $\sigma=-5.47559\dots$  so this will be sufficient. However, in this paper we give a different and complete proof of \eqref{forsigma}  for $\sigma<\frac12$.

This double integral appears in one of the papers of Riemann Nachlass in Göttingen. Usually we are able to interpret what Riemann does in any of these papers, but in this case, they have not aided me in understanding what Siegel says. 

Not only for the above considerations, but rather by the transformations that Siegel writes between equations (36) and (54) of his paper, we wonder from time to time if the Riemann papers that Siegel consulted returned to their place in the Nachlass. There are several papers in the actual Nachlass with the Riemann-Siegel expansion, but Siegel does not appear to be following any of them in his paper. Nevertheless, Siegel says: \emph{Just how strong Riemann's analytic technique was, it is made particularly clear by his derivation and manipulation of the asymptotic series for $\zeta(s)$.}

\subsection*{16. Density theorems for Riemann's auxiliary function \texorpdfstring{\cite{A174}}{[21]}}
Iwaniec and Kowalski  (\cite{IK}*{p.~2}) say that density theorems about the zeros of L-functions are likely a waste of time studying the empty set. Rather, these density theorems  probably cannot be better because they are really theorems about the zeros of $\Rzeta(s)$, and this function has many zeros spread over the entire critical strip.

\begin{figure}[H]
\begin{center}
\includegraphics[width=0.99\hsize]{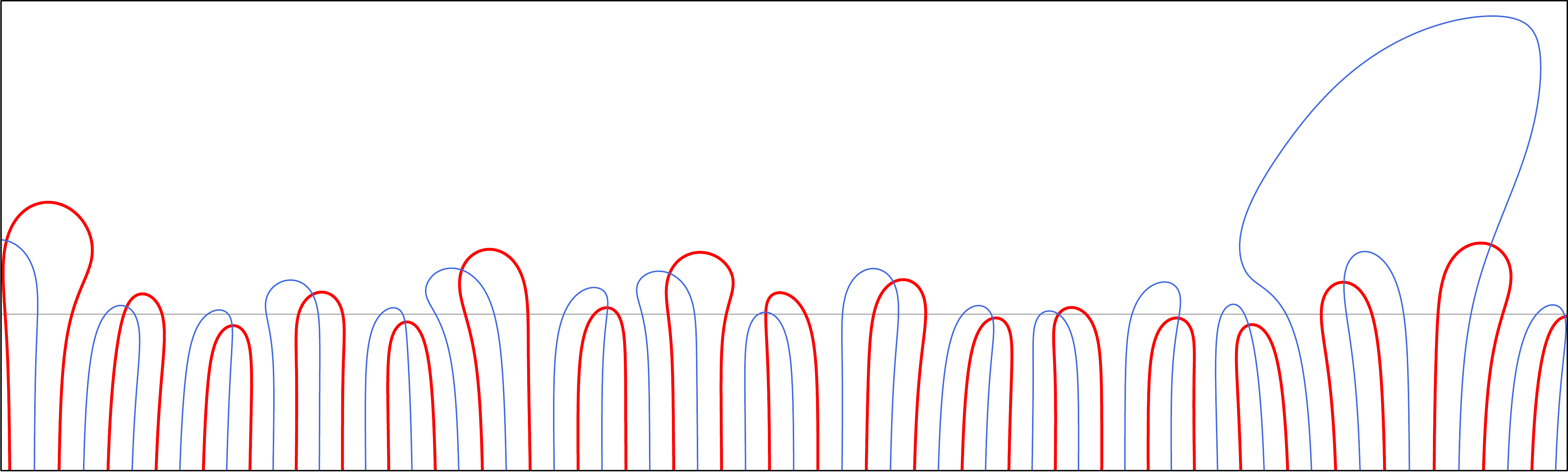}
\caption{x-ray of $e^{i\vartheta(t)}\Rzeta(\frac12+it)$ in $(200\,040,200\,060)\times(-2,4)$}
\label{blueredfigure}
\end{center}
\end{figure}

Since 
\[Z(t)=e^{i\vartheta(t)}\zeta(\tfrac12+it)=2\Re\{e^{i\vartheta(t)}\Rzeta (\tfrac12+it)\}.\]
$\zeta(s)$ have a zero on the critical line at each point $s=\frac12+it$  where $e^{i\vartheta(t)}\Rzeta (\tfrac12+it)$  is purely imaginary.  This is the same as saying that
$\zeta(s)$ have a zero at a point $s$ in the critical line just in case $\pi^{-s/2}\Gamma(s/2)\Rzeta(s)$ is purely imaginary.  
Siegel connected the zeros of $\Rzeta(s)$ with zeros on the critical line of $\zeta(s)$, showing that $\zeta(s)$ have at least $2N_2(T)$ zeros on the critical line with ordinates between $0$ and $T$, where $N_2(T)$ is the number of zeros of $\Rzeta(s)$ to the left of the critical line and with ordinates between $0$ and $T$. To understand the relationship between the zeros, consider the x-ray of 
$e^{i\vartheta(t)}\Rzeta(\frac12+it)$ in Figure \ref{blueredfigure}. The horizontal line is the real axis for $t$.

Each zero of $\Rzeta(s)$ to the left of the critical line gives a zero above the real axis in Figure \ref{blueredfigure}. The function $Z(t)$ vanishes at the points where the imaginary lines (blue lines) in the figure cut the real axis. 
All imaginary lines (in blue) start at the bottom of the figure \ref{blueredfigure} so that the imaginary line passing through a zero of $e^{i\vartheta(t)}\Rzeta(\frac12+it)$
above the real axis cuts the real axis at two points, producing two zeros of $Z(t)$. 

In principle, it would be possible that the imaginary line through a zero below the real axis in the figure does not cut the real axis. But we see in the examples in the figure that they are sufficiently close to give also two zeros of $Z(t)$.

\begin{figure}[H]
\begin{center}
\includegraphics[width=0.99\hsize]{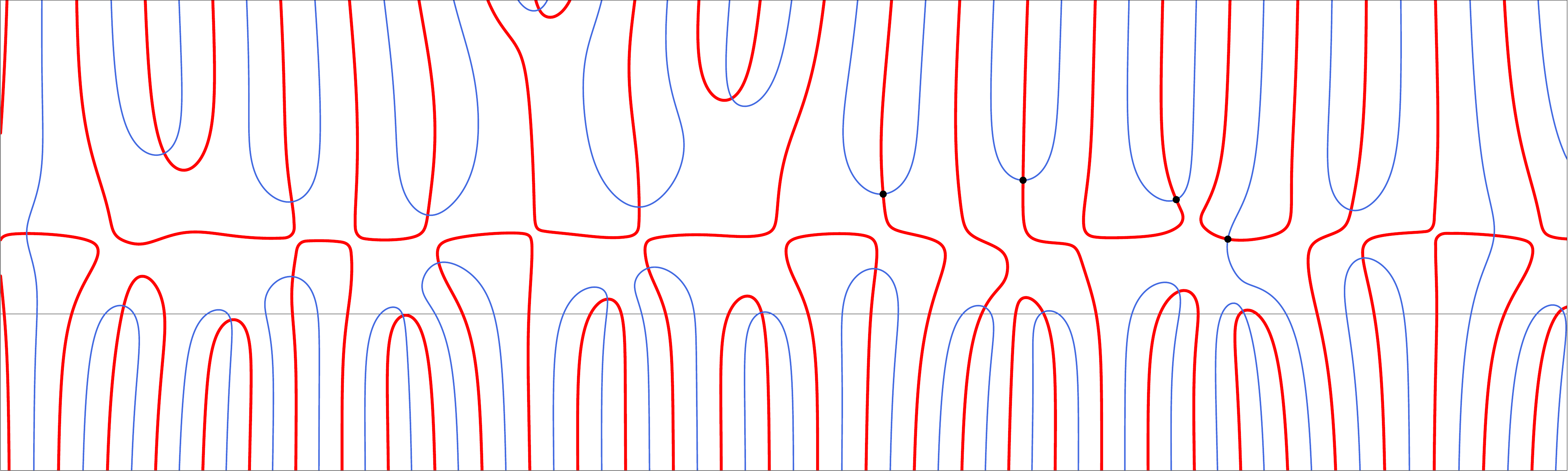}
\caption{x-ray of $e^{i\vartheta(t)}G(t)$ in $(200\,040,200\,060)\times(-2,4)$}
\label{blueredfigureG}
\end{center}
\end{figure}

Since $\Rzeta(s)$ can be approximated by  Dirichlet polynomials, the methods used to prove the density theorems for $\zeta(s)$ can be applied to $\Rzeta(s)$, and this is the objective of this paper. Our main Theorem  proves that for $1/2<\alpha\le1$, we have $N(\alpha,T)\ll T^{\frac32-\alpha}\log^3T$, where $N(\alpha,T)$ denotes the number of zeros of $\Rzeta(s)$ in $[\alpha,1]\times[0,T]$ counted with multiplicity.  The proof is standard, but we do not see any advantage in using a mollifier here. Perhaps that is a point where my result can be easily improved.  

These zeros below the real line in our figure correspond to the zeros of $\Rzeta(s)$ to the right of the critical line. In this paper, we show that most of the zeros of $\Rzeta(s)$ to the right of the critical line are  near this line. These density theorems  make it almost inevitable that each of these zeros produce its two zeros for $\zeta(s)$. In this way, the $\frac{T}{4\pi}\log\frac{T}{2\pi}-\frac{T}{4\pi}$ zeros of $\Rzeta(s)$ can generate the $\frac{T}{2\pi}\log\frac{T}{2\pi}-\frac{T}{2\pi}$ zeros  predicted by the Riemann hypothesis.

In this paper, we follow the Theorems in Iwaniec \cite{Iw},
adapting them to  $\Rzeta(s)$. We give proofs of these Theorems, giving explicit constants. 

\subsection*{17. Left density theorems for the auxiliary function \texorpdfstring{\cite{A191}}{[22]}}
Since $\zeta(s)-\Rzeta(s)=\chi(s)\overline{\Rzeta(1-\overline s)}$ the zeros to the left of the critical line are symmetrical with respect to this line of the zeros of $\zeta(s)-\Rzeta(s)$. We may approximate this last function by a Dirichlet polynomial to the right of the critical line. This may be a way to bound the number of these zeros. This is one of the papers in which my lack of technique is most exposed. We have not been able to obtain anything reasonable following this line. We mention this in case anyone can get something.

\subsection*{18. On the approximation of the zeta function by Dirichlet polynomials \texorpdfstring{\cite{A192}}{[23]}} 
For its use in \cite{A191}, we needed an explicit version of Theorem 4.11 of Titchmarsh. There are many papers dealing with this question, but we found a real chaos in the literature on the range of validity of the various results. Therefore, we prove that for $s=\sigma+it$ with $\sigma\ge0$ and $0<t\le x$, we have $\zeta(s)=\sum_{n\le x}n^{-s}+x^{1-s}/(s-1)+\frac{29}{14}\Theta x^{-\sigma}$, where $\Theta$ is a complex number with $|\Theta|\le1$. This improves Theorem 4.11 of Titchmarsh, making explicit a result from Iwaniec in \cite{Iw} to whom we dedicate this work. In the end, we do not arrive at anything good in \cite{A191}. 

\section{Representation and approximation to Riemann-Siegel function \texorpdfstring{$Z(t)$}{Z(t)}}

There are many functions $f(s)$ satisfying $Z(t)=2\Re\{e^{i\vartheta(t)}f(\frac12+it)\}$.
For example, in \cite{A62} with van de Lune we proved
\[-\vartheta'(t)Z(t)=\Re\{e^{i\vartheta(t)}\zeta'(\tfrac12+it)\}.\]
Each of these formulas gives an opportunity to show that many zeros are on the critical line. Many times we get an approximation to $f(\frac12+it)$ giving us an approximation to $Z(t)$. The general objective here is to obtain a real function $u(t)$ such that $Z(t)=u(t)+\orden(t^{-\varepsilon})$, and such that $u$ has $\frac{T}{2\pi}\log\frac{T}{2\pi}-\frac{T}{2\pi}+\orden(T)$ zeros on $[0,T]$. 

\subsection*{19. An expression for Riemann Siegel function \texorpdfstring{\cite{A75}}{[24]}} Lavrik \cite{Lav} consider for any $\tau\in\C$ with $\Re\tau>0$ the function 
\begin{equation}\label{E:seriesrep}
\Lzeta(\tau,s)=-\frac{(\pi\tau)^{s/2}}{s\;\Gamma(\frac{s}{2})}+\sum_{n=1}^\infty \frac{\Gamma(\frac{s}{2},\pi n^2\tau)}{\Gamma(\frac{s}{2})}\frac{1}{n^s}.
\end{equation}
This function satisfies 
\begin{equation}\label{E:InRiemann}
\pi^{-s/2}\Gamma(\tfrac{s}{2})\zeta(s)=\pi^{-s/2}\Gamma(\tfrac{s}{2})\Lzeta(\tau,s)+\pi^{-(1-s)/2}\Gamma(\tfrac{1-s}{2})\Lzeta(1/\tau,1-s).
\end{equation}
We consider specially the case $\tau=1$ that was considered by Riemann in \cite{R}. We define $\Lzeta(s)=\Lzeta(1,s)$, and we have $Z(t)=2\Re\{e^{i\vartheta(t)}\Lzeta(\frac12+it)\}$. The non-trivial zeros of $\Lzeta(s)$ are contained in the halfplane
$\Re s>11.25$ and are easy to compute. 

\subsection*{20. Approximate formula for \texorpdfstring{$Z(t)$}{Z(t)}\texorpdfstring{\cite{A60}}{[25]}} 
\begin{figure}[H]
\begin{center}
\includegraphics[width=0.99\hsize]{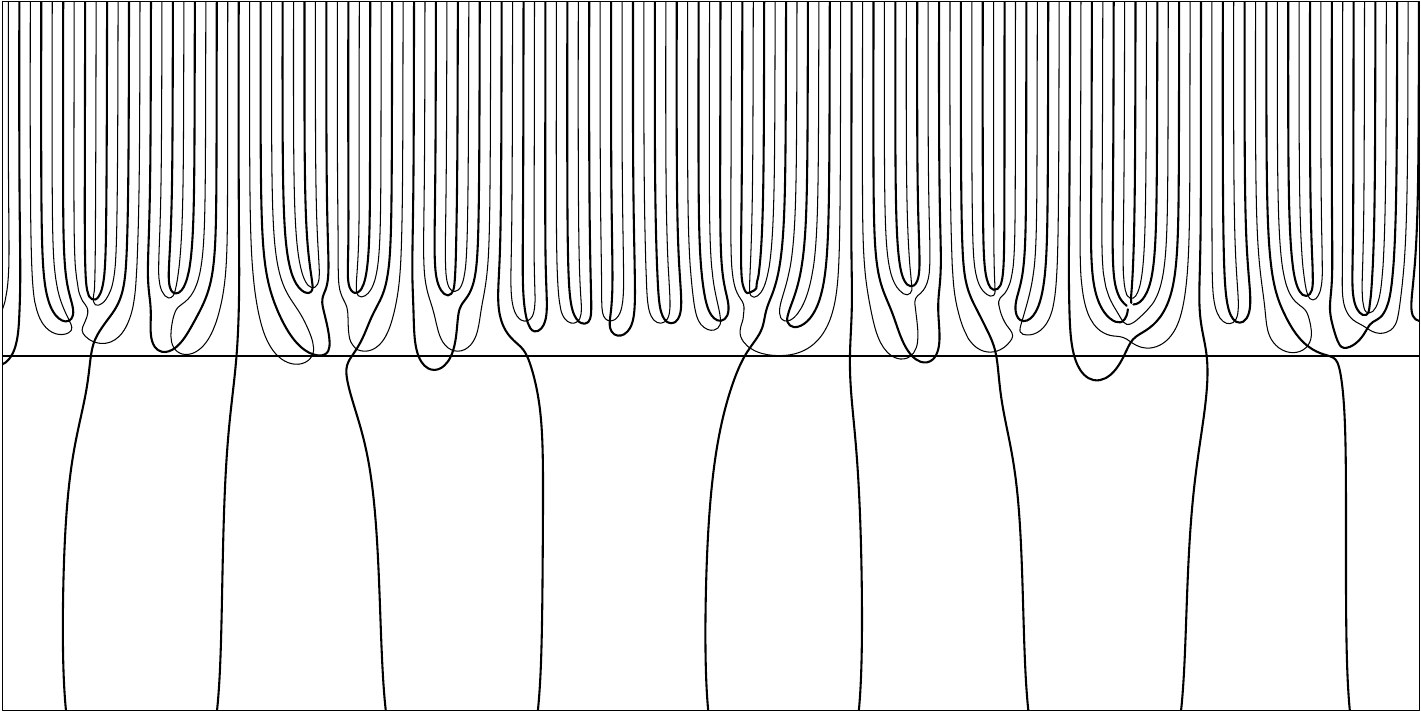}
\caption{Plot of  $G(t)$ in $(1000,1040)\times(-10,10)$.}
\label{secondplot3}
\end{center}
\end{figure}

The function
\begin{equation}\label{defU}
U(t)=\frac{1}{2\pi i}\int\limits_{-i\sigma-\infty}^{-i\sigma+\infty}
e^{i\vartheta(t+x)-i\vartheta(t)}\zeta\bigl(\tfrac12+i(t+x)\bigr)\frac{\pi}{2\sinh\frac{\pi x}{2}}\,dx,
\end{equation}
where $0<\sigma<\frac12$, satisfies $Z(t)=2\Re\{e^{i\vartheta(t)}U(t)\}$. Approximating this function, we get 
\begin{equation}\label{F:G}
G(t)=\sum_{n=1}^\infty \frac{1}{n^{\frac12+it}}\frac{t}{2\pi n^2+t}
\end{equation} satisfies $Z(t)=2\Re\{e^{i\vartheta(t)}G(t)\}+\Orden(t^{-\frac56+\varepsilon})$. 
The function $G(t)$ extends to a meromorphic function on $\C$ with poles at $(2k-\frac12)i$ and $-2\pi n^2$. There are many interesting things about $G(t)$ for example
\begin{equation}
G(t)=\frac{t}{2\pi}\Bigl(\zeta(2+\tfrac12+it)-\sum_{n=1}^\infty \frac{1}{n^{2+\frac12+it}}
\frac{t}{2\pi n^2+t}\Bigr).
\end{equation}
So that the values of $Z(t)$ appear to depend only on things happening at $\sigma>2$.
Also for small $t$ we observe in Figure \ref{secondplot3} that the continuous  $\arg G(t)$ is increasing very slowly. Giving the zeros we want for $Z(t)$.

\subsection*{21. Integral Representation for Riemann-Siegel \texorpdfstring{$Z(t)$}{Z(t)}  function \texorpdfstring{\cite{A162}}{[26]}} We apply the solution of the Dirichlet problem in the strip $-4\le\sigma\le 5$ for adequate boundary values to get an integral representation of $Z(t)$
\begin{equation}\label{E:intform2}
F(t)=\int_{-\infty}^\infty 
\frac{h(x)\zeta(4+ix)}{7\cosh\pi\frac{x-t}{7}}\,dx, \qquad Z(t)=\frac{\Re F(t)}{(\frac14+t^2)^{\frac12}(\frac{25}{4}+t^2)^{\frac12}},
\end{equation}
where $h(x)$ is a simple function. After that we get the estimate
\begin{equation}
Z(t)=\Bigl(\frac{t}{2\pi}\Bigr)^{\frac74}\Re\{e^{i\theta(t)}H(t)\}+\Orden(t^{-3/4}).,
\end{equation}
with 
\[H(t)=\int_{-\infty}^\infty\Bigl(\frac{t}{2\pi}\Bigr)^{ix/2}\frac{\zeta(4+it+ix)}{7\cosh(\pi x/7)}\,dx=\Bigl(\frac{t}{2\pi}\Bigr)^{-7/4}\sum_{n=1}^\infty \frac{1}{n^{\frac12+it}}\frac{2}{1+(\frac{t}{2\pi n^2})^{-7/2}}.\]
We arrive at an equation similar to \eqref{F:G}. The behaviour of this $H$ is similar to that of $G$. 

\section{Some integrals for \texorpdfstring{$\zeta(s)$}{zeta(s)}, \texorpdfstring{$\Rzeta(s)$}{R(s)} and similar functions}

\subsection*{22. An entire function defined by Riemann  \texorpdfstring{\cite{A152}}{[27]}}
This is a representation integral for $\zeta(s)$ that we find in one of the papers of Riemann's Nachlass in Göttingen\footnote{We have a promise to communicate to Göttingen any papers I obtained from the papers of Riemann, we expect to do it before publishing this.}.
Riemann defines a bounded function on $\R$ whose Fourier transform vanishes at the real values with $\zeta(\frac12+i\gamma)=0$. First, Riemann defines an entire function
\begin{equation}\label{e1}
\chi(z):=\int_L\frac{e^{izw}}{e^{2\pi w^2}-1}\,dw,
\end{equation}
where $L$ is a line parallel to the real axis at a height $0<b<2^{-1/2}$. It has the power series expansion
\begin{equation}\label{e3}
\chi(z)=\frac{z}{2}+\frac{1}{\sqrt{2}}\sum_{n=0}^\infty \frac{\zeta(n+\frac12)}{n!}\Bigl(-\frac{z^2}{8\pi}
\Bigr)^n.
\end{equation}
and satisfies the functional equation $\chi(-z)=\chi(z)-z$. And is connected with $\zeta(s)$ by 
\begin{equation}\label{e8}
\int_0^\infty \chi(x)x^{2s-1}\,dx=\Gamma(2s)\zeta(s+\tfrac12)\cos\tfrac{\pi(2s+1)}{4},
\qquad \Re(s)>0.
\end{equation}

\subsection*{23. An integral representation of \texorpdfstring{$\Rzeta(s)$}{R(s)} due to Gabcke \texorpdfstring{\cite{A187}}{[28]}}

Gabcke \cite{G} proved
\[\Rzeta(s)=2^{s/2}\pi^{s/2}e^{\pi i(s-1)/4}\int_{-\frac12\searrow\frac12}
\frac{e^{-\pi i u^2/2+\pi i u}}{2i\cos\pi u}U(s-\tfrac12,\sqrt{2\pi}e^{\pi i/4}u)\,du\]
where $U(\nu,z)$ is the usual parabolic cylinder function \cite{AS}. Here, we give a new proof and put it in the form
\begin{equation}\label{withhermite}
\Rzeta(s)=-2^s \pi^{s/2}e^{\pi i s/4}\int_{-\infty}^\infty \frac{e^{-\pi x^2}H_{-s}(x\sqrt{\pi})}{1+e^{-2\pi\omega x}}\,dx.
\end{equation}
where $H_\nu(z)$ is the generalization of Hermite polynomials \cite{L}*{S.~10.2}.

\subsection*{24. Integral representation of Riemann auxiliary function \texorpdfstring{\cite{A188}}{[29]}}

We transform \eqref{withhermite}  into 
\begin{equation}\label{FirstInt}
\Rzeta(s)=-\frac{2^s \pi^{s}e^{\pi i s/4}}{\Gamma(s)}\int_0^\infty
y^{s-1}\frac{1-e^{-\pi y^2+\pi \omega y}}{1-e^{2\pi \omega y}}\,dy.
\end{equation}
What is remarkable here is that the function of the integrand $f(y)=\frac{1-e^{-\pi y^2+\pi \omega y}}{1-e^{2\pi \omega y}}$ appearing in \eqref{FirstInt} is entire. Zeros in the denominator are also zeros in the numerator. For $y\to+\infty$ we have $|f(y)|\sim e^{-\pi y\sqrt{2}}$. When we move the integration line, no residues appear. This may be useful in bounding $\Rzeta(s)$ without introducing the zeta sum. The function $f(y)$ has a nice x-ray,  shown in the paper.

\subsection*{25. Levinson functions \texorpdfstring{\cite{A197}}{[30]}}  Levinson \cite{L2} find entire functions of $s$, $f(s,\tau)$ for each $\tau\ne0$ on the closed hiperplane $\Im\tau\ge 0$ (my notations in this paper), such that 
\[h(s)\zeta(s)=h(s) f(s,\tau)+h(1-s)\overline{f}(1-s,1/\overline\tau), \qquad h(s)=\pi^{-s/2}\Gamma(\tfrac{s}{2}).\]
This implies that the function 
\[\Rzeta_\tau(s):=\Rzeta(s,\tau):=\frac{f(s,\tau)+f(s,1/\overline\tau)}{2},\]
satisfies 
\begin{equation}\label{mainexamples}
Z(t)=2\Re\{e^{i\vartheta(t)}\Rzeta_\tau(\tfrac12+it)\}.
\end{equation}
We also have $\Rzeta_{-1}(s)=\Rzeta(s)$ and $\Rzeta_i(s)=\Lzeta(s)$ considered in \cite{A75}. The limit cases when $\tau$ is a non-null rational number are especially interesting. The definition given by Levinson is by integrating a series. In this paper we show that for rational $\tau$, this series can be simplified to a trigonometric function as in the case of $\Rzeta$. To achieve this we prove some nice summation formulas. Similar to Gaussian sums. For example, we show:
\begin{equation}
\begin{aligned}
\Rzeta_{-4/3}(s)=&\frac{1}{2\sqrt{3}}\int_{0\swarrow1}\frac{x^{-s}e^{\frac{4\pi i}{3}x^2}}{e^{\pi i x}-e^{-\pi i x}}\Bigl(2e^{-\frac{\pi i}{6}}\cos\tfrac{\pi x}{3}+i\cos(\pi x)\Bigr)\,dx\\
&+\frac{1}{4}\int_{0\swarrow1}\frac{x^{-s}e^{\frac{3\pi i}{4}x^2}}{e^{\pi i x}-e^{-\pi i x}}\Bigl(e^{\frac{\pi i}{4}}+2\cos\tfrac{\pi x}{2}+e^{-\frac{3\pi i}{4}}\cos(\pi x)\Bigr)
\end{aligned}
\end{equation}

\begin{equation}
\begin{aligned}
\Rzeta_{-2}(s)=\frac12\int_{0\swarrow1}&\frac{x^{-s}e^{2\pi ix^2}\cos\pi x}{e^{\pi i x}-e^{-\pi i x}}\,dx\\&+\frac{1}{2\sqrt{2}}\int_{0\swarrow1}\frac{x^{-s}e^{\frac{\pi i}{2}x^2}}{e^{\pi i x}-e^{-\pi i x}}\Bigl(e^{\frac{\pi i}{4}}+e^{-\frac{\pi i}{4}}\cos(\pi x)\Bigr)\,dx.
\end{aligned}
\end{equation}
This must be joined with \eqref{mainexamples}. In the paper we give general similar expressions for any rational number not equal to $0$.

The above expressions hint that we may transform the approximate functional equation in Theorem 4.13 of Titchmarsh into asymptotic expansions. 

\section{A naive integral}

\subsection*{26. A naive integral \texorpdfstring{\cite{A123}}{[33]}}
In \cite{A166} two real functions $g(x,t)$ and $f(x,t)$ are defined, so that the Riemann-Siegel $Z$ function is given as 
\[Z(t)=\Re\Bigl\{\frac{u(t)e^{\frac{\pi i}{8}}}{\frac12+it}\int_0^\infty g(x,t)e^{i f(x,t)}\,dt\Bigr\},\]
where $u(t)$ is a real function of order $t^{-1/4}$ when $t\to+\infty$. The function $g(x,t)$ is indefinitely differentiable and tends to $0$ as well as all its derivatives when $x\to0^+$ or $x\to+\infty$. Since, furthermore, for $t\to+\infty$ the function $f(x,t)$ tends to $+\infty$  we may expect that the integral depends essentially on the behavior of $g(x,t)$ at the extremes. 

As Polya in an analogous situation \cite{P} we consider the substitution of $\psi(x)$ by a simpler  similar function.
A simple function with this behavior is 
\[\psi_0(x):=2\pi(1+\tfrac{1}{4}x^{-5/2})e^{-\pi x-\frac{\pi}{4x}}.\]
Therefore, we define $J_0(t)$ replacing in the definition  of $J(t)$ the function $\psi(x)$ by  the simpler $\psi_0(x)$.
\begin{equation}
J_0(t)=2\pi\int_0^\infty (1+\tfrac{1}{4}x^{-\frac52})e^{-\pi x-\frac{\pi}{4x}}(1-ix)^{\frac12(\frac12+it)}\,dx.
\end{equation}
The resulting $Z_0(t)$ disappoints us
\[Z_0(t)\asymp 
\Re\Bigl\{\frac{2}{\sqrt{\pi}}\exp\Bigl\{i\Bigl(\frac{t}{2}\log\frac{t}{2\pi}-\frac{t}{2}-
\frac{\pi}{8}\Bigr)\Bigr\}+\frac{2}{(2\pi t)^{1/4}}\exp\Bigl(\pi i\sqrt{\frac{t}{2\pi}}\;\Bigr)\Bigr\},\quad t\to+\infty.\]
However, the integral $J_0(t)$ is interesting as a technical challenge. And still we have the possibility to get a better result improving $\psi_0(x)$.  

This paper needs a strong revision, but we will upload it to arXiv without this revision. It is an interesting integral and can be a good challenge for those interested in the asymptotic development of integrals.

\section{Some other related results}

\subsection*{27. Simple bounds for the auxiliary function of Riemann \texorpdfstring{\cite{A92}}{[31]}}
Given $s\ne1$ in $\C$ and $\varepsilon>0$, in \cite{A88} we want to compute $z$ such that $|\zeta(s)-z|<\varepsilon$. To this end we need to compute many products $ab$ with a prescribed error. To achieve this we need some rough bounds of the quantities $a$ and $b$. This is the origin of this paper. Its object is to give bounds of $\zeta(s)$, $\vartheta(s)$, $\Rzeta(s)$, $Z(t)$, for use in the numerical computations. It is good to have reliable, not necessarily very good, but easy to compute bounds. The objective were not to publish this, but only serve me for reference. Since we  are not going to publish these papers and we have many references, we prefer to upload it as it currently stands.

\subsection*{28. Explicit van der Corput estimates \texorpdfstring{\cite{A181}}{[32]}}

In \cite{A173} it was required to bound a zeta sum. We required explicit bounds, so the book of Titchmarsh was not sufficient. We noticed that many papers overlooked a result by Landau. This made us upload \cite{A182} to arXiv, which is nothing more than a translation of the paper in German by Landau, where he proves what is the best constant for an inequality which some pretended to improve upon. 

Then, we searched for theorems to bound the zeta sum that can be turned explicit. The obvious election was van der Corput's $d$-th estimates. We searched for the original papers and found, besides the usual references,  three papers of van der Corput giving explicit estimates. It was surprising that no one referred to these papers. It is a very elaborate construction, so we started studying the papers. We found an error in the van der Corput reasoning. This explains why nobody refers to these papers.  We found that the error can be rectified fairly easily. This paper is the exposition of van der Corput reasoning on these papers with adequate correction.  We prove
\begin{theorem} \label{T:dthtest}
Let $X$, and $Y\in\R$ be such that $\lfloor Y\rfloor>d$ where $d\ge3$ is a natural number. 
Let $f\colon(X,X+Y]\to\R$ be a real function with continuous derivatives up to the order $d$. Assume that $0<\lambda\le f^{(d)}(x)\le\Lambda$ for $X<x\le X+Y$.  Denote by $D=2^d$. Then 
\begin{equation}\label{E:mainmain}
\Bigl|\frac{1}{Y}\sum_{X<n\le X+Y}e(f(n))\Bigr|\le
\max\Bigl\{A_d\Bigl(\frac{\Lambda}{\lambda Y}\Bigr)^{2/D}, B_d\Bigl(\frac{\Lambda^2}{\lambda}\Bigr)^{1/(D-2)},C_d(\lambda Y^d)^{-2/D}\Bigr\},
\end{equation}
where 
$A_d$, $B_d$, and $C_d$ are explicit constants. They depend on $d$ but for $d\ge2$ for example $A_d< 7.5$, $B_d<5.8$ and $C_d<10.9$.
\end{theorem}
We apply this theorem to zeta sums, giving the best choice of $d$ in each case. Also, we prove that our Theorem implies Titchmarsh's Theorem 5.13.

\subsection*{29. Levinson-Conrey method following the exposition of Iwaniec \texorpdfstring{\cite{A194}}{[34]}}
Contains a detailed exposition of the Levinson-Conrey method. This is the result of my study of the book by Iwaniec \cite{Iw}. Contains several changes to the book and improvements, some of them due to Iwaniec.  I appreciate his answers to my many questions. My aim was to apply this method directly to the $\Rzeta(s)$ function, but I have not reached any useful result. This is almost a booklet in which I have worked much,but it is as much my work as it is Iwaniec's.

\subsection*{30. Programs for computing the zeros \texorpdfstring{\cite{A180}}{[35]}}

It is difficult to get all the zeros $\rho=\beta+i\gamma$ of $\Rzeta(s)$  with $0<\gamma\le T$. 
Here, we explain the procedure we have followed. But the paper is incomplete, and consequently, I will not upload it for the moment.

\subsection*{31. X-ray's of holomorphic functions}

Still, only a project. It will give definitions, some interesting properties, and how to read an x-ray. It will contain some history, for example that Gauss published the first x-ray, as an illustration of his first proof of the fundamental theorem of algebra. The complex plane and the x-rays appeared at the same time.

\section{Some Open Problems}\label{S:prob}

An outside observer will see here nothing but methods already known, when a new idea is surely needed to advance the proof of Riemann's last theorem. Actually, this is just the result of new ideas not getting me anywhere, as, for example, in \cite{A191}, but I always thought that everything we knew about $\Rzeta(s)$ could be useful to get to the difficult proof that Riemann claimed to have. I have not even finished expressing the already known methods; for example, I could write a paper explaining the relation of $\Rzeta(s)$ and the Lindelöf hypothesis. The Lindelöf hypothesis is equivalent to saying that $\Rzeta(\frac12+it)\ll t^\varepsilon$ for all $\varepsilon>0$. Therefore, Riemann's efforts to bound $\Rzeta(\frac12+it)$ in \eqref{E:doubleInt}, could lead him to prove Lindelöf's hypothesis.

I also want to point out here the difference between the density theorems for $\Rzeta(s)$ and $\zeta(s)$. Putting the zeros of $\zeta(s)$ close to the critical line does not seem to move in the direction of proving Riemann's last theorem, but doing the same with $\Rzeta(s)$ is different.  Every zero of $\Rzeta(s)$ close to the line is likely to place two zeros of $\zeta(s)$ right on the line. 

Here we collect some open problems related to $\Rzeta(s)$. 

\begin{problem}
Can we get some useful bound from the expression of $|\Rzeta(\sigma+it)|^2$ as a double integral.
\end{problem}
I mention here that the very similar double integrals of Riemann in his Nachlass are not equal to the ones given by Siegel.

\begin{problem}
Let $\Rzeta(\rho)=0$ with $\Im(\gamma)>0$. Is it true that $\Re(\beta)<1$?
\end{problem}

This is proved in \cite{A173} for $\Im(\gamma)>3.93\times10^{65}$, and is true for $\Im(\gamma)<215946$.

Another important problem is the left limit of the zeros. Siegel \cite{Siegel} claim to have proved it. But we have not been able to prove it, so we consider it an open problem.
\begin{problem} Let $\varepsilon>0$. Is it true that there is no zero of $\Rzeta(s)$ in the region  defined by  $1-\sigma>t^\varepsilon$ and $t>t_0$, for some $t_0>0$?
\end{problem}

\begin{problem}\label{wellzeros}
Let $N_{\Rzeta}(T)$ be the number of zeros of $\Rzeta(s)$ with $0<\Im(\rho)\le T$, and $N_{\Rzeta}(\beta<1/2,T)$, the number of those zeros with $\Re(\rho)<1/2$, both counted with multiplicities. Does
\[\lim_{T\to+\infty}\frac{N_{\Rzeta}(\beta<1/2,T)}{N_{\Rzeta}(T)}\]
exist? Is its value equal to 1?
\end{problem}
An affirmative answer to the second question in Problem \ref{wellzeros} implies the Riemann last theorem. 

One possibility of proving the Riemann last theorem is that it would be a very general property. 

\begin{problem}
Let $f$ be a holomorphic function for $\Re s>0$ with $|f(s)-1|<1/2$ for $\sigma>2$ and $t>t_0$. Assume also that 
\[\int_0^T|f(\tfrac12+it)|^2\,dt=\Orden( T\log T),\]
and that  there is some constant $a$ such that  $f(s)=\Orden(t^a)$ for $0<\sigma<2$ and $t\to+\infty$. 

Is it true that the real function $\Re\{e^{i\vartheta(t)}f(\frac12+it)\}$ has $\sim\frac{T}{2\pi}\log\frac{T}{2\pi}$ zeros on the segment $(0,T)$?
\end{problem}
Here we have to be open about the exact conditions for $f$.

\section*{Afterword}

These twenty-odd articles represent my unsuccessful attempts to prove Riemann's last theorem. As long as I did not achieve my goal, I did not publish my progress. This has always been my way of doing mathematics. But at 76 years of age, I realised that I did not want these small advances to be definitively lost with my death. So a few years ago I started to write them more carefully. With the idea of uploading them to arXiv. 

I don't achieve my goal, but I bring some knowledge about the function $\Rzeta(s)$, and I pose several problems that I think are interesting (see section \ref{S:prob}).
The main reason for publishing them is that I still think that there is a real way to prove the last Riemann theorem by means of the $\Rzeta(s)$ function. Sure some mathematician will be able to achieve what I have not been able to. 

I have reviewed these papers as if I were the referee. In such a large volume I'm sure some error has escaped me, but I hope there are no serious problems. I can continue reviewing, but there will be no end to it.  There is one article that goes nowhere \cite{A191} and another one \cite{A123} that needs revision and I'm not going to do that because it's just a side issue.   That's why I'm deciding to upload them now. 

I would be appreciate any comments that could improve a possible new version.

\end{document}